\documentclass[reqno, 11pt]{amsart}
\usepackage{fullpage}  
\usepackage{bbold,bbm,color,calc,graphicx,amsfonts,amsthm,amscd,epsfig,psfrag,amsmath,amssymb,enumerate,dsfont}
\usepackage{caption}
\usepackage{subcaption}

\usepackage{pdfsync}
\usepackage{stmaryrd} 
\usepackage[numeric,initials,nobysame]{amsrefs}
\usepackage{marginnote}
\usepackage[top=2.5cm, bottom=2.5cm, outer=2.5cm, inner=2.5cm, heightrounded, marginparwidth=1.9cm, marginparsep=0.3cm]{geometry}

\usepackage[colorlinks=true, pdfstartview=FitV, linkcolor=blue,citecolor=blue, urlcolor=blue,pagebackref=false]{hyperref}
\usepackage[showlabels,sections,floats,textmath,displaymath]{}
\usepackage{booktabs}
\usepackage{appendix}
\usepackage[percent]{overpic}
\newcommand{\beq}{\begin{equation}}
\newcommand{\eeq}{\end{equation}}
\newcommand{\ie}{\hbox{\it i.e.\ }}

\renewcommand{\restriction}{\mathord{\upharpoonright}}

\reversemarginpar
\newlength\fullwidth
\numberwithin{equation}{section}

\DeclareMathSymbol{\leqslant}{\mathalpha}{AMSa}{"36} 
\DeclareMathSymbol{\geqslant}{\mathalpha}{AMSa}{"3E} 
\DeclareMathSymbol{\eset}{\mathalpha}{AMSb}{"3F}     
\renewcommand{\leq}{\;\leqslant\;}                   
\renewcommand{\geq}{\;\geqslant\;}                   

\DeclareSymbolFont{bbold}{U}{bbold}{m}{n}
\DeclareSymbolFontAlphabet{\mathbbold}{bbold}

\renewcommand{\b}{\beta}

\newcommand{\var}{\operatorname{Var}}

\newcommand{\tmix}{T_{\rm mix}}
\newcommand{\trel}{T_{\rm rel}}

\newcommand{\D}{\Delta}

\renewcommand{\b}{\beta}
\renewcommand{\l}{\lambda}
\renewcommand{\L}{\Lambda}

\renewcommand{\l}{\lambda}
\renewcommand{\a}{\alpha}
\renewcommand{\d}{\delta}
\renewcommand{\t}{\tau}

\newcommand{\g}{\gamma}

\newcommand{\e}{\varepsilon}

\renewcommand{\o}{\omega}
\renewcommand{\O}{\Omega}

\newcommand{\gap}{{\rm gap}}
\newcommand{\tc}{\thinspace |\thinspace}

\newtheorem{theorem}{Theorem}[section]
\newtheorem{pseudo-theorem}{Pseudo-Theorem}

\newtheorem{lemma}[theorem]{Lemma}
\newtheorem{proposition}[theorem]{Proposition}
\newtheorem{corollary}[theorem]{Corollary}
\newtheorem{remark}[theorem]{Remark}

\newtheorem{definition}[theorem]{Definition}

\newcommand{\N}{\mathbb N}

\newcommand{\cA}{\ensuremath{\mathcal A}}
\newcommand{\cB}{\ensuremath{\mathcal B}}
\newcommand{\cC}{\ensuremath{\mathcal C}}

\newcommand{\cF}{\ensuremath{\mathcal F}}
\newcommand{\cG}{\ensuremath{\mathcal G}}

\newcommand{\cL}{\ensuremath{\mathcal L}}

\newcommand{\cN}{\ensuremath{\mathcal N}}

\newcommand{\cP}{\ensuremath{\mathcal P}}

\newcommand{\cT}{\ensuremath{\mathcal T}}

\newcommand{\bbE}{{\ensuremath{\mathbb E}} }

\newcommand{\bbN}{{\ensuremath{\mathbb N}} }

\newcommand{\bbP}{{\ensuremath{\mathbb P}} }

\newcommand{\bbR}{{\ensuremath{\mathbb R}} }

\newcommand{\bbZ}{{\ensuremath{\mathbb Z}} }

\let\a=\alpha \let\b=\beta   \let\d=\delta  \let\e=\varepsilon
 \let\g=\gamma \let\h=\eta      \let\l=\lambda
      \let\o=\omega      
  \let\s=\sigma \let\t=\tau   
  
\let\D=\Delta   \let\G=\Gamma  \let\L=\Lambda 
\let\O=\Omega      

\renewcommand{\le}{\leq}

\usepackage[normalem]{ulem}

\definecolor{light}{gray}{.9}

\title[Mixing time and local ergodicity of the East-like process in
$\bbZ^d$]{Mixing time and local exponential ergodicity of the East-like process in $\bbZ^d$} 
\author[P. Chleboun]{P. Chleboun}
 \address{P. Chleboun. Mathematics Institute,
Zeeman Building, University of Warwick, Coventry CV4 7AL.} \email{p.i.chleboun@warwick.ac.uk}

\author[A. Faggionato]{A. Faggionato}
\address{A. Faggionato. Dipartimento di  Matematica, Universit\`a  La  Sapienza, P.le Aldo Moro  2, 00185  Roma, Italy. }\email{faggiona@mat.uniroma1.it}

 \author[F. Martinelli]{F. Martinelli}
 \address{F. Martinelli. Dipartimento di  Matematica e Fisica,  Universit\`a  Roma Tre, Largo S.L.Murialdo 00146, Roma, Italy. } \email{martin@mat.uniroma3.it}

\begin{document}
\begin{abstract}
  The East process, a well known reversible linear chain of spins,
  represents the prototype of a general class of interacting particle
  systems with constraints modeling
  the dynamics of real glasses. In this paper we consider a
  generalization of the East process living in the $d$-dimensional
  lattice and we establish new progresses on the \emph{out-of-equilibrium}
  behavior. In particular we prove a form of (local) exponential
  ergodicity when the initial distribution is far from the stationary
  one and we prove that the mixing time in a finite box grows linearly
  in the side of the box.
\end{abstract}
\maketitle
\section{Introduction}
Kinetically constrained spin models (KCMs) are interacting
$0$-$1$ particle systems, on general graphs, which evolve with a simple Glauber dynamics
described as follows. At every site $x$ the system tries to update the occupancy variable (or spin) at
$x$ to the value $1$ or $0$ with probability $p$ and $q=1-p$\footnote{In
  the physical applications $q\approx e^{-c \beta}$ at low
  temperature, where $\b$ is the inverse-temperature and $c$ is a
  constant.},
respectively. However the
update at $x$ is accepted only if the \emph{current} local
configuration satisfies a certain constraint, hence the models are ``kinetically
constrained''. It is always assumed that the constraint at site $x$ does
not depend on the spin at $x$ and therefore the product
Bernoulli($p$) probability measure $\pi$ is the reversible measure.  Constraints
may require, for example, that a certain \emph{number} of the
neighbouring spins are in state $0$, or more restrictively, that certain
\emph{preassigned} neighbouring spins are in state $0$ (e.g. the
children of $x$ when the underlying graph is a rooted tree). 

The main interest in the physical literature for KCMs (see e.g. 
\cites{Ritort,GarrahanSollichToninelli} for a review) stems from the
fact  that they display many key dynamical
features of real glassy materials: ergodicity breaking
transition at some critical value $q_c$, huge relaxation time for $q$
close to $q_c$, dynamic heterogeneity (\ie non-trivial spatio-temporal fluctuations of the
local relaxation to equilibrium) and aging, just to mention a
few. Mathematically, despite their simple definition, KCMs pose
very challenging and interesting problems because of the hardness of
the constraint, with ramifications towards bootstrap percolation
problems \cite{Spiral}, combinatorics \cites{CDG,Valiant:2004cb},
coalescence processes \cites{FMRT-cmp,FMRT} and  random walks on
upper triangular matrices \cite{Peres-Sly}. Some of the mathematical tools
developed for the analysis of the relaxation process of KCMs
\cite{CMRT} have proved to
be quite powerful also in other contexts such as card shuffling problems
\cite{Bhatnagar:2007tr} and random evolution of surfaces \cite{PietroCaputo:2012vl}.

Among the KCMs, the most studied model is the \emph{East process}, a
one-dimensional spin system that was introduced in the physics
literature by J\"{a}ckle and Eisinger~\cite{JACKLE} (cf. also \cites{SE1,SE2}
and \cite{East-Rassegna} for a recent mathematical review). In this case the
base graph is $\bbZ$ (or finite connected subsets of it) and the
constraint at $x\in \bbZ$ requires that the vertex $x-1$ is empty (\ie
its associated spin is $0$).

It is the properties of the East process before and towards reaching equilibrium
which are of interest, with the standard gauges for the speed of
convergence to stationarity given by the relaxation time $\trel$
($\equiv$ inverse spectral-gap) and the total-variation mixing time $\tmix$
on a finite interval $\L=\{0,\ldots,L\}$ (cf. e.g. \cite{Levin2008}), where we fix $\eta(0)=0$ for ergodicity.

That 
$\trel=O(1)$\footnote{We recall that $f = O(g)$ means that there exists a
  constant $C > 0$ such that $|f|\le C |g|$ and that
  $f=\O(g)\Leftrightarrow g=O(f)$. } in $L$ for any
$q$ small enough was first proved in \cite{Aldous} and, later on, for
all $q\in (0,1]$ in 
\cite{CMRT} by different methods.
Subsequently the analysis of the relaxation time in the physical relevant setting $q\searrow 0$
(corresponding to the low temperature limit) was developed to a
high level of precision in \cites{CFM,CFM-JSTAT} where the relevant questions of
dynamical heterogeneities and time scale separation have been rigorously settled.
The fact that the relaxation time is $O(1)$ in $L$
implies, in particular, that $\tmix\sim L$ (cf. Theorem \ref{th:mix}). It is then natural to ask whether
the finite volume East process exhibits the \emph{cutoff phenomenon}
(coined by Aldous and Diaconis~\cite{AD86}; see
also~\cites{Aldous,DiSh,Diaconis0} and the references therein): over a negligible
period of time (the cutoff window) the distance from equilibrium drops
from near $1$ to near $0$. As is easily seen, the cutoff problem is strongly linked to
the following \emph{front
progression} problem for the East process.
For a initial configuration $\eta$ with $\sup\{x:\eta_x=0\}<\infty$,
call this rightmost 0 its \emph{front} and denote by $X(\eta)$ its position. It is easy to verify
that at any later time the process starting from $\eta$ will also have
a front $X(\eta_t)$ and the key step to prove cutoff is a detailed
analysis of the asymptotic law of
$X(\eta_t)$ as $t\to \infty$. In \cite{Blondel} a kind
of shape theorem was proved: as $t\to\infty$  the law of the process behind the front 
converges to an invariant measure $\nu$ and $\frac1{t} X(\eta_t) \to
v_\infty$ in probability for a suitable constant $v_\infty>0$. As a consequence $\tmix \sim
L/v_\infty$. To prove cutoff one has to go beyond
the law of large numbers and to control the fluctuations of the
front around the mean value $v_\infty t$. In \cite{East-cutoff} it was
proved that the latters obey a CLT, and cutoff with the optimal window
$\sim \sqrt{L}$ follows.

In several interesting contributions (cf. \cites{Gar1, Gar2,Garrahan2003} and references therein)
a natural generalization of the East dynamics
to higher dimensions $d>1$, in the sequel referred to as the \emph{East-like}
process, appears to play a key role in \emph{realistic}
models of glass formers. In $d=2$ the East-like process evolves 
similarly  to the East process but now the kinetic
constraint requires that the South \emph{or} West neighbor of the
updating vertex contains \emph{at least} one vacancy (\ie a zero spin). In general
$\eta(x)$ can flip if $\eta(x-e)=0$ for some $e$ in the
canonical basis of $\bbZ^d$. 

In \cites{CFM2,CFM2-EPL} the authors thoroughly analyzed the East-like
process with emphasis on its low temperature (small $q$)
behavior. Among the main results it was proved that the process on $\bbZ^d$ is ergodic with a finite relaxation time
$\trel(\bbZ^d)$ and that, as $q\searrow 0$, $\trel(\bbZ^d)\asymp
\trel(\bbZ)^{1/d}$, correcting in this way some heuristic conjectures
which appeared in the physical literature.  In finite volume with ergodic boundary conditions (see Section
\ref{sec:fin-vol} below) the asymptotics of the relaxation time as
$q\searrow 0$ was also computed quite precisely and it was shown to depend very strongly on
the choice of the boundary conditions. Finally the asymptotics as
$q\searrow 0$ of the
\emph{persistence times} (see Section \ref{sec:persistence-time}
below) was also established but only up to a spatial scale
$O(1/q^{1/d})$. 

Other natural questions concerning the 
genuine \emph{out-of-equilibrium} behavior of the East-like process,
for example the convergence as $t\to \infty$ of the particle density
to the equilibrium value $p$ 
for a general class of initial laws or the asymptotics as $L\to
\infty$ of the mixing time $\tmix$ in
a box of side $L$ remained open. 

In this paper we fill this gap in the spirit of similar results for
the East process \cite{CMST}. In Section 4 (cf. Theorem \ref{th:1} and
Corollary \ref{cor:1}), using key results from
Section 3, we
prove a form of \emph{local exponential ergodicity} together with the above mentioned convergence of the particle
density. In the same section we also prove an exponential tail for the law of the
persistence time of a given vertex (cf. Theorem \ref{th:2}), a key
step for our final main result, namely  that 
the mixing time in a finite box  grows linearly in the side of the
box (cf. Theorem \ref{th:mix}).

\section{The Model}
\subsection{Notation}
For any $x=(x_1,\dots,x_d)\in \bbZ^d$ we denote its $\ell^1$-norm by
$\|x\|_1=\sum_{i=1}^d|x_i|$. For any $x\in \bbZ^d$ we define the
positive quadrant $\bbZ^d_{x,+}$ at $x$ as the set $\{y\in \bbZ^d:\
y_i\ge x_i,\ i=1,\dots,d\}$.  When $x$ is the orgin we will simply
write $\bbZ^d_+$. We will also let $\bbZ^d_{x,\uparrow}\triangleq
\bbZ^d_{x,+}\setminus \{x\}$ and $\bbZ^d_{x,\downarrow}\triangleq
\{y\in \bbZ^d:\, x\in \bbZ^d_{y,\uparrow}\}$. The interpretation of
these two sets is as follows:  $\bbZ^d_{x,\uparrow}$ is
the set of vertices other than $x$ which are \emph{influenced} by the
spin at $x$, while $\bbZ^d_{x,\downarrow}$ are those vertices other
than $x$ which \emph{influence} the spin at $x$.

We denote by $\cB=\{e_1, e_2 , \dots ,e_d\}$ the canonical basis of
$\bbZ^d$ and, given a set $X \subset \bbZ^d$, we define its
East-boundary by
$$ 
\partial_E X\triangleq \{ y \in \bbZ^d \setminus X\,:\, x+e_i \in X \text{ for
  some } e_i\}\,.
$$
Given $\L\subset \bbZ^d$, we will denote by $\O_\L$ the product space
$ \{0,1\}^\L$ endowed with the product topology. If $\L = \bbZ^d$ we
simply write $\O$. In the sequel we will refer to the vertices of $\L$
where a given configuration $\eta\in \O_\L$ is equal to one (zero) as
the \emph{particles} (\emph{vacancies}) of $\eta$.  If $V\subset \L$
and $\eta\in \O_\L$ we will write $\eta\restriction_{V}$ for the
restriction of $\eta$ to $V$. In particular we will simply write
$\eta\restriction_V=1$ to mean that $\eta(x)=1$ for all $x\in V$.  
Finally, for any $\L\subset \bbZ^d$, a configuration
$\s\in \O_{\partial_E \L}$ will be referred to as a \emph{boundary
  condition}.

\subsection{Finite volume process and boundary conditions}
\label{sec:fin-vol}
Given a region $\L \subset \bbZ^d$ and a configuration
$\s\in \O_{\partial_E \L}$, we define the \emph{constraint} at site $x\in
\L$ with boundary condition $\s$, in the sequel denoted by $c_x^{\L,
  \s} (\eta)$, as the indicator function of the event in $\O_\L$ that there exists
$e\in \cB$ such that $(\s\cdot\h)(x-e)=0$,   
where $\s\cdot\h \in \Omega_{\L \cup \partial_E \L}$ is the configuration equal
to $\s$ on $\partial_E \L$ and $\h$ on $\L$.  Then the East-like
process with parameter $q \in (0,1)$ and boundary configuration $\s$
is the continuous time Markov chain with state space $\O_\L$ and
infinitesimal generator
\begin{align}
  \label{eq:gen}
  \cL^\s _\L f(\eta) &= \sum _{x \in \L} c_x^{\L, \s}(\eta)\, \bigl[
  \eta_x q+ (1-\eta_x)p \bigr] \cdot \bigl[ f(\eta^x) -f(\eta) \bigr]
  \nonumber\\
  &= \sum _{x \in \L} c_x^{\L, \s}(\eta) \bigl[\pi_x (f)- f\bigr]
  (\eta),
\end{align}
where $p= 1-q$, $\eta^x$ is the configuration in $\O_\L$ obtained
from $\eta$ by flipping its value at $x$, and $\pi_x$ is the
Bernoulli($p$) probability measure on the spin at $x$.

We call the boundary condition $\sigma$ \emph{ergodic} if the process
given by \eqref{eq:gen} is ergodic.  For example if $\L= \prod_{i=1}^d
[a_i,b_i] $ then a boundary condition $\s$ is ergodic if and only if
$\s(a-e) = 0$ for some $e \in \cB$.  If $\s$ is such that by removing
some single vacancy in $\s$ one obtains a non-ergodic boundary
condition then $\s$ is said to be \emph{minimal}.  On the other hand
if $\s$ is identically equal to zero we call it a \emph{maximal}
boundary condition.

Since the local constraint $c_x^{\L, \s}(\eta)$ does not depend on
$\eta_x$, it is simple to check
that the East-like process is reversible w.r.t. the product
Bernoulli($p$) probability measure $\pi_\L=\prod_{x\in \L}\pi_x$ on $\O_\L$. In
what follows we will usually write $\pi$ instead of $\pi_\L$ whenever
no confusion can arise. Moreover,
since our results are uniform in the choice of the ergodic boundary
condition, we will usually drop the superscript $\s$ in our notation
and we will write $\bbP^{\L}_\eta(\cdot)$ and
$\bbE^{\L}_\eta(\cdot)$ for the law and the associated expectation of
the process with initial condition $\eta$. The only exception is when $\Lambda = \bbZ_+^d$ and the boundary condition $\s$ is minimal, \ie
$\s(-e) = 0$ for some $e\in\cB$ and $\s(x) = 1$ otherwise. In that
case we will write $\bbP^{\rm min}_{\eta}(\cdot)$ and $\bbE^{\rm min}_{\eta}(\cdot)$.

\subsection{The infinite volume East-like process}
We also define  the East process on the entire lattice $\bbZ^d$ as follows. 
Let
$
c_x (\eta) \triangleq \mathbbold{1}{(\exists\, e \in \cB \text{ such that } \eta(x-e)=0)},
$
be the constraint at $x\,$\footnote{We will adopt the notation
  $\mathbbold{1}{(A)}$ for the indicator of the event $A$}.
Then the East-like process on $\bbZ^d$ is the continuous time Markov
process with state space $\O$, with reversible  measure given by the
product Bernoulli($p$) probability measure $\pi=\prod_{x\in \bbZ^d}\pi_x$ and infinitesimal generator $\cL$ whose
action on functions depending on finitely many spins is given by  
\begin{align}
 \cL f(\eta) &= \sum _{x \in \bbZ^d} c_x(\eta)\,  \bigl[ \eta_x q+ (1-\eta_x)p \bigr] \cdot
\bigl[ f(\eta^x) -f(\eta) \bigr]\nonumber\\
\label{eq:generator}&=\sum _{x \in \bbZ^d} c_x(\eta)\bigl[\pi_x (f)- f\bigr] (\eta)\,.
 \end{align} 
We will denote by $\bbP_\eta(\cdot)$ and
$\bbE_\eta(\cdot)$ the law and the associated expectation of
the process started from $\eta$. In general, if the law of the initial
configuration is $\nu$, we will write $\bbP_\nu(\cdot),\, \bbE_\nu(\cdot)$
\subsection{Spectral Gap}  
In the sequel, the spectral
gap of a reversible Markov generator $\cA$ will be denoted by $\gap(\cA)$.
It is defined as infimum over all non-constant functions $f$ in the
domain of $\cA$ of the ratio between the Dirichlet form of $f$ and its
variance w.r.t. the reversible probability measure. In the finite
dimensional case, if $\cA$ is irreducible then $\gap(\cA)>0$. A positive spectral gap implies, in particular, 
 that the reversible measure is
mixing for the semigroup generated by $\cA$ and that the variance
contracts exponentially fast:
\begin{align}
  \var(e^{t\cA}f)\leq e^{-2 t\,\gap(\cA)}\var(f), \quad
  \forall\, f \in L^2\,.
\end{align}
We now recall some properties of the spectral gap for the East-like
process which will be important in the following (for more details see
 \cite{CFM2,CFM2-EPL}).
The infinite volume East-like process has a positive spectral gap for
all $p \in (0,1)$. If $\L\subset \bbZ^d$ then $\gap(\cL_{\L}^\s)$ is bounded from below by a positive constant uniformly
in choice of the boundary condition $\s$ among the ergodic  ones. Moreover, if 
$\L=\prod_{i=1}^d[1,\ell_i]$, then $\gap(\cL_{\L}^\s)$ is decreasing
in $\{\ell_i\}_{i=1}^d$ \emph{and}  in the values of the boundary spins
$\{\s(x)\}_{x\in \partial_E\L}$.

\subsection{Graphical construction}
\label{sec:graph}
We recall a graphical construction of the East–like process on $\bbZ^d$, which will be very useful in the sequel. A similar construction, with slight modifications, holds also in the finite volume case. To each $x\in\bbZ^d$, we associate a rate one Poisson process and,
independently, a family of independent Bernoulli$(p)$ random variables
$\{s_{x,\ell} : \ell \in \N\}$ ($\bbN:= \{0,1,\dots\}$). The occurrences of the Poisson process
associated to $x$ will be denoted by $ \{t_{x,\ell} : \ell \in \N\}$,
labelled in increasing order.  We assume independence as $x$ varies in
$\bbZ^d$. Sometimes in the sequel we will refer to the collection
$\{t_{x,\ell},s_{x,\ell}\}_{\ell\in \bbN}$ as the \emph{clock rings
  and coin tosses} associated to the vertex $x$. We write $(\Theta,
\bbP)$ for the probability space on which the above objects are
defined.  Notice that, $\bbP$-almost surely, all the occurrences
$\{t_{x,\ell} : \ell \in \N,\, x\in \bbZ^d\}$ are different (this
property will be often used below without further mention).

Given a probability measure $\nu$ on $\O$ we consider the product
probability measure $\bbP_\nu\triangleq \nu \otimes \bbP$ on the product space
$\O \otimes \Theta$. Then on $(\O \otimes \Theta,\bbP_\nu)$ we can
define a c\`adl\`ag Markov process $(\eta_t) _{t \geq 0}$ which is exactly the
East-like process on $\bbZ^d$ given by \eqref{eq:generator} with initial distribution $\nu$, as follows. 
 The initial configuration
$\eta_0$ associated to the element $(\eta, \vartheta)\in \O \otimes
\Theta$ is given by $\eta$.  At each time $t=t_{x,\ell}(\vartheta)$
the site $x$ queries the state of its own constraint $c_x\bigl(
\eta_{t-}\bigr)$.  If and only if the constraint is satisfied (\ie
$c_x\bigl(\eta_{t-}\bigr) = 1$), then $t $ is called a \emph{legal ring}
and the configuration $\eta_{t}$ is obtained from $\eta_{t-}$ by
resetting its value at site $x$ to the value of the corresponding
Bernoulli variable $s_{x,\ell}(\vartheta)$. Using the Harris's percolation
argument \cite{Harris72} the above definition is well posed for $\bbP_\nu$
a.e. $(\eta, \vartheta)$. When $\nu= \d_\eta$ we simply write
$\bbP_\eta$ instead of $\bbP_\nu$\footnote{Although we have defined $\mathbb{P}_\eta , \bbP_\nu$ also for the law of the East--like process, it will be clear from the context when the notation is referred to the above graphical construction.}.

Finally, for $\L\subset \bbZ^d$ and $t>0$, we let $\cF_\L$
($\cF_{\L,t}$) be the $\s$-algebra generated by all
clock rings and coin tosses associated to vertices in $\L$ (all
the clock rings and coin tosses up to
time $t$).
\section{Local stationarity} In this section we prove two results
showing a kind of
local stationarity of the reversible measure $\pi$. In the first case
the region where we want to prove stationarity of  $\pi$ is a
non-random subset $\L$ of $\bbZ^d$. In the second case 
the set $\L$ will be random and
determined by the dynamics itself in the graphical construction (the definition uses clock rings and coin tosses).  

\begin{proposition}
\label{lem:0}
Let $\L$ be a finite subset of $\bbZ^d$ and assume that the initial distribution $\nu$ of the East-like process in
$\bbZ^d$ is the product of its marginals on $\O_\L,\O_{\L^c}$ and that the marginal on $\O_\L$ coincides
with $\pi_\L$. Then, for any $t>0$ and any $\s\in\O_\L$, 
\[
\bbP_\nu(\eta_t\restriction_\L=\s\tc \cF_{\L^c,t})=\pi(\s).
\]
\end{proposition}
\begin{proof}
Let $x^*\in \L$ be such that $\bbZ^d_{x^*,\uparrow}\cap \L=\emptyset$, let $\L^*\triangleq \L\setminus \{x^*\}$
and let $\cF_t^{*}\triangleq
\cF_{\bbZ^d\setminus \{x^*\},t}$. Clearly such a vertex always
exists. Then
\begin{align*}
&\bbP_\nu(\eta_t\restriction_\L=\s\tc
\cF_{\L^c,t})=\int d\nu(\eta)\bbE_\eta\left(\mathbbold{1}{(\eta_t\restriction_{\L^*}=\s\restriction_{\L^*})}\bbP_\eta\bigl(\eta_t(x^*)=\s(x^*)\tc \cF_t^{*}\bigr)\tc \cF_{\L^c,t}\right) \\
&=\int d\nu(\eta)\bbE_\eta\Bigl(\mathbbold{1}{(\eta_t\restriction_{\L^*}=\s\restriction_{\L^*})}\sum_{\xi\in \{0,1\}}\pi(\eta(x^*)=\xi)\bbP_{\eta,\xi}\bigl(\eta_t(x^*)=\s(x^*)\tc \cF_t^{*}\bigr) \tc \cF_{\L^c,t}\Bigr),
  \end{align*}
  where the notation $\bbP_{\eta,\xi}(\cdot)$ indicates that the
  initial configuration is equal to $\eta$ outside $x^*$ and equal to
  $\xi$ at $x^*$.  Above we used the definition of $x^*$ to guarantee
  that, once the initial condition $\eta$ is given, the event
    $\{\eta_t\restriction_{\L^*}=\s\restriction_{\L^*}\}$
  is measurable w.r.t. $\cF_t^{^*}$. Also, by the same reason, the
  event does not depend on the value $\eta(x^*)$. Finally, we have used the
  assumption on $\nu$ to perform a partial average over the initial
  value of $\eta(x^*)$.

We now claim that
\[
\sum_{\xi\in
  \{0,1\}}\pi(\eta(x^*)=\xi)\bbP_{\eta,\xi}(\eta_t(x^*)=\s(x^*)\tc
\cF_t^*)=\pi(\s(x^*)).
\]
To prove the claim we condition on the event $\cA_{t}$ that there has been
at least one legal ring at
$x^*$ before time $t$. Notice that also $\cA_{t}$ does not depend on the
initial value $\eta(x^*)$. Thus
\begin{align*}
\sum_{\xi\in
  \{0,1\}}&\pi(\eta(x^*)=\xi)\bbP_{\eta,\xi}(\eta_t(x^*)=\s(x^*)\tc
\cF_t^*)\\
&=
\sum_{\xi\in
  \{0,1\}}\pi(\eta(x^*)=\xi)\bbP_{\eta,\xi}(\eta_t(x^*)=\s(x^*)\tc
\cA_{t},\cF_t^*) \bbP_\eta(\cA_{t}\tc \cF_t^*)\\
&+ \sum_{\xi\in
  \{0,1\}}\pi(\eta(x^*)=\xi)\bbP_{\eta,\xi}(\eta_t(x^*)=\s(x^*)\tc
\cA^c_{t},\cF_t^*) \bbP_\eta(\cA^c_{t}\tc \cF_t^*)
\end{align*}
We now observe that, for any $\eta$, 
\[
(i)\quad \bbP_{\eta,\xi}(\eta_t(x^*)=\s(x^*)\tc
\cA_{t},\cF_t^*) = \pi(\s(x^*)),
\]
because in this case $\eta_t(x^*)$ takes the value of the last coin toss at $x^*$
before $t$, and 
\[
(ii)\quad \sum_{\xi\in
  \{0,1\}}\pi(\eta(x^*)=\xi)\bbP_{\eta,\xi}(\eta_t(x^*)=\s(x^*)\tc
\cA^c_{t},\cF_t^*)=\pi(\s(x^*)),
\]
because $\eta_t(x^*)=\xi$ on the event $\cA_{t}^c$. Hence the claim.

In conclusion
\begin{align*}
\bbP_\nu(\eta_t\restriction_\L=\s\tc \cF_{\L^c,t})&= 
\bbP_\nu(\eta_t\restriction_{\L^*}=\s\restriction_{\L^*}\tc
\cF_{\L^c,t})\ \pi(\s(x^*))\\
&=\bbP_\nu(\eta_t\restriction_{\L^*}=\s\restriction_{\L^*}\tc
\cF_{(\L^*)^c,t}) \ \pi(\s(x^*)).
\end{align*}
Since $\pi$ is a product measure, the term 
$\bbP_\nu(\eta_t\restriction_{\L^*}=\s\restriction_{\L^*}\tc
\cF_{(\L^*)^c,t})  $ has the required form on the
reduced set $\L^*$ and the proof follows by iteration.
\end{proof}

\begin{remark} In the sequel we will use the above proposition but
  only in the ``easy'' case in which the set $\L$ has empty
  intersection with $\bbZ^d_{x,\downarrow}$ for any
  $x \in \partial_E \L$ (e.g.  $\L$ is given by a box,
  $\bbZ^d_+$ or some $\bbZ^d_{y,\uparrow}$).  In this case, under the
  same assumptions of Proposition \ref{lem:0}, the dynamics in
  $\partial _E \L$ never queries the state of the dynamics in $\L$.
  Moreover, calling $t_1 <t_2 < \cdots < t_n $ the times at which a
  spin of $\partial _E \L$ flips in the time window $[0,t]$ and
  setting $t_0\equiv 0$, $t_{n+1} \equiv t$, for any $i=0,1, \dots, n$
  in the time interval $[t_i, t_{i+1})$ the projection on $\L$ of the
  East--like process equals a.s.  the East--like process on $\L$ with
  fixed boundary condition $\eta \restriction_{\partial_E \L} (t_i)$.
\end{remark}

We now extend Proposition \ref{lem:0} to a case in which the
set $\L$ is itself random. Given a
realization of the clock  rings and coin tosses at all vertices and an initial configuration $\eta$ such that
$\eta(0)=0$, let $\t_0=0, z^{(0)}=0$  and define
\begin{align*}
\t_{k+1}&=\inf\{s>\t_k:\ \text{at time $s$ there is a legal ring at
  $z^{(k)}$}\}\,,\\
z^{(k+1)}&=\min\{x\in \partial_E\{z^{(k)}\}:\   \eta_{\t_{k+1}^-}(x) =0\} \,,
\end{align*}
where the minimum is  taken w.r.t. the lexicographic order.  Notice
that $\eta_{t}(z^{(k)})=0$  for any $ k$, $t \in [\t_k ,\t_{k+1})$  a.s.  We will refer to
this special vacancy as the \emph{distinguished zero}
(cf. \cite{Aldous} and \cite{East-Rassegna} for an analogous
definition for the one dimensional East process). 
  \begin{definition}
\label{def:1} 
We
define the \emph{trace} of the distinguished zero up to time $t$ as the set
\[  \G_t\triangleq \left\{z^{(0)}, z^{(1)}, \dots, z^{(\cN_t-1)}\right\},\]
where $\cN_t\triangleq \max \{ k \geq 0\,:\, \t_k \leq t\}$ and we use the convention that $\G_t = \emptyset $ if $\cN_t=0$.
\end{definition}
\noindent
The above definition is well posed $\bbP_\eta$--a.s.  since
$\lim_{k \to \infty } \t_k=+\infty$ $\bbP_\eta$--a.s.  
Let also
$\cG_0^{(1)}$ be the $\s$-algebra containing all information ``below'' the origin, up to the
first legal ring $\t_1$. Formally $\cG_0^{(1)}$ is generated by all
events $F$ such that
$F\cap \{\t_1\le s\}\in \cF_{\bbZ^d_{0,\downarrow},s}$ for any
$s \geq 0$. Note that the event $\{z^{(1)}=z\}$ belongs to
$\cG_0^{(1)}$.

The above construction is of course valid for any initial vacancy of
the initial configuration $\eta$. If the distinguished zero is
initially at $x$ then we will simply add a subscript ``$x$'' to the above
notation.

We will now construct recursively
the $\s$-algebras $\cG_x^{(n)}$, $x\in \bbZ^d$, containing all the information
on the first $n\ge 2$ steps of a distinguished zero initially at $x$ and also on the clock rings and
coin tosses ``below'' the successive positions of the distinguished zero and
between consecutive jumps. 
In order to proceed more formally we first need the following definition.
\begin{definition}[Time shift in $\Theta$]
Given $s>0$ together with $\o\triangleq
\{t_{y,j},s_{y,j}\}_{y\in \bbZ^d,j\in \bbN}\in \Theta$, 
let
\[
\theta_s\o\triangleq\{t_{y,j+\nu_{y,s}},s_{y,j+\nu_{y,s}}\}_{y\in \bbZ^d,j\in \bbN},
\] 
where $\nu_{y,s}\triangleq \min \{j:\ t_{y,j}\ge s\}$ (recall that $\bbN=\{0,1,\dots\}$). In other words
the first ring and coin toss at $y$ for $\theta_s\o$ are the first
ring and coin toss at $y$ \emph{after}
$s$ and so on.  
\end{definition}
We then define recursively the family $\{\cG_x^{(n)}\}_{n\ge 2}$ as
follows. $\cG_x^{(n)}$ is the $\s$-algebra generated by all events of
  the form 
  \begin{equation}
    \label{eq:7}
F^{(n)}=F^{(1)}\cap\{z^{(1)}=z\}\cap \{\theta_{\t_1}\o\in
F^{(n-1)}\},\quad z\in \partial_E\{x\},
  \end{equation}
where $F^{(1)}\in \cG_x^{(1)}$ and $F^{(n-1)}\in \cG_z^{(n-1)}$.

We are finally ready to state our main result on the law of $\eta_t\restriction_{\G_t}$.
\begin{proposition}
\label{prop:miracoli}
For all $n\ge 1$, all $\eta$ such that
$\eta(0)=0$ and all $t>0$ the conditional distribution of
$\eta_t\restriction_{\G_t}$ given $\cG_0^{(n)}$ and $\{\cN_t=n\}$ coincides with the reversible measure
$\pi$. 
\end{proposition}
\begin{proof}
Fix
$n\ge 1$, $\s=(\s_0,\dots,\s_{n-1})\in \{0,1\}^n$ and
$z\in \partial_{E}\{0\}$. Let also $F^{(n)}\in \cG_0^{(n)}$ be as in \eqref{eq:7}.
It is enough to prove that\footnote{ We write $\bbP(A;B)$ instead of $\bbP(A\cap B)$
  for shortness.}
\begin{gather}
\bbP_\eta(\cN_t=n;\ F^{(n)};\ \eta_t\restriction_{\G_t}=\s)=\pi(\s) \bbP_\eta(\cN_t=n;\ F^{(n)}).
  \label{eq:AA1}
\end{gather}
We will exploit \eqref{eq:7} and use induction. Let $\hat \eta_{\t_1}$ be the configuration in
$\bbZ^d\setminus\{0\}$ equal to
$\eta_{\t_1}$ at all vertices different from the origin. Let also
$\bbP_{\hat\eta_{\t_1},\pi}(\cdot)$ be the law of the East-like process
with initial condition equal to $\hat\eta_{\t_1}$ outside the origin and
sampled from $\pi$ at the origin. The strong Markov
property\footnote{Here we appeal to the strong Markov property of the
underlying Poisson point process given by the clock rings.}
w.r.t. the Markov time $\t_1$ together with \eqref{eq:7} imply that
\begin{align}
\label{eq:AA2}
\bbP_\eta(\cN_t=n;\ &F^{(n)};\ \eta_t\restriction_{\G_t}=\s)
\nonumber\\
=\bbE_\eta\Bigl[&\mathbbold{1}{(F^{(1)}\cap\{z^{(1)}=z\})}\mathbbold{1}{(\t_1<t)}\
                  \times\\
&\bbP_{\hat\eta_{\t_1},\pi}\left(\cN_{z,t-\t_1}=n-1;\ F^{(n-1)};\
  \eta_{t-\t_1}\restriction_{\G_{z,t-\t_1}}=\s';\ \eta_{t-\t_1}(0)=\s_0\right)\Bigr],
\nonumber
\end{align}
where $\s'=(\s_1,\dots,\s_{n-1})$ and we adopt the convention that the
event $\{\eta_{t-\t_1}\restriction_{\G_{z,t-\t_1}}=\s'\}$ is the
sure event if $\G_{z,t-\t_1}=\emptyset$. 

Observe now that, given $\hat\eta_{\t_1}$, for all $s>0$ the events
$\{\cN_{z,s}=n-1\},\, F^{(n-1)}$ and $\{\eta_{s}\restriction_{\G_{z,s}}=\s'\}$ are
measurable w.r.t. $\cF_{\bbZ^d_{0,\downarrow},s}$ while the event
$\{\eta_{s}(0)=\s_0\}$ is measurable
w.r.t. $\cF_{\bbZ^d_{0,\downarrow}\cup\{0\},s}$. Therefore, using
Proposition \ref{lem:0} applied to the set $\L=\{0\}$,
\begin{gather}
\bbP_{\hat\eta_{\t_1},\pi}\left(\cN_{z,s}=n-1;\ F^{(n-1)};\
  \eta_{s}\restriction_{\G_{z,s}}=\s';\ \eta_{s}(0)=\s_0\right)\nonumber\\
=\bbE_{\hat\eta_{\t},\pi}\left[\mathbbold{1}{(\cN_{z,s}=0)}\mathbbold{1}{(F^{(n-1)})}\mathbbold{1}{(\eta_{s}\restriction_{\G_{z,s}}=\s')}\bbP_{\hat\eta_{\t},\pi}(\eta_s(0)=\s_0\tc
 \cF_{\{0\}^c,s} )\right]\nonumber\\
=\pi(\s_0)\bbP_{\hat\eta_{\t},\pi}(\cN_{x,s}=0;\ F^{(n-1)};\ \eta_{s}\restriction_{\G_{z,s}}=\s').
\label{eq:AA3}
\end{gather}
If we now apply \eqref{eq:AA3} to \eqref{eq:AA2} we get 
\begin{gather*}
\bbP_\eta(\cN_t=n;\ F^{(n)};\ \eta_t\restriction_{\G_t}=\s)\\ 
=
\pi(\s_0) \bbE_\eta\left[\mathbbold{1}{(F^{(1)}\cap\{z^{(1)}=z\})}\mathbbold{1}{(\t_1<t)}
\bbP_{\hat\eta_{\t},\pi}(\cN_{z,t-\t_1}=n-1;\ F^{(n-1)};\ \eta_{t-\t_1}\restriction_{\G_{z,t-\t_1}}=\s')\right]
\end{gather*}
and \eqref{eq:AA1} follows for $n=1$. 
If $n\ge 2$ and we inductively assume \eqref{eq:AA1} for
$n-1$ to write
\[
\bbP_{\hat\eta_{\t},\pi}(\cN_{z,t-\t_1}=n-1;\ F^{(n-1)};\ \eta_{t-\t_1}\restriction_{\G_{z,t-\t_1}}=\s')=\pi(\s')
\bbP_{\hat\eta_{\t},\pi}(\cN_{z,t-\t_1}=n-1;\ F^{(n-1)}),
\]  
we get
\begin{gather*}
\bbP_\eta(\cN_t=n;\ F^{(n)};\ \eta_t\restriction_{\G_t}=\s)\\
=\pi(\s) \bbE_\eta\left[\,\mathbbold{1}{(F^{(1)}\cap\{z^{(1)}=z\})}\mathbbold{1}{(\t_1<t)}
\bbP_{\hat\eta_{\t},\pi}(\cN_{z,t-\t_1}=n-1;\ F^{(n-1)})\,\right]\\
=\pi(\s)\bbP_\eta\left(\cN_t=n;\ F^{(n)}\right),
\end{gather*}
\ie \eqref{eq:AA1} for $n$.
 \end{proof}

\section{Out-of-equilibrium results}
In this section we begin by proving two results (cf. Lemma \ref{lem:0bis} and
Corollary \ref{cor:1bis} below) showing that, given an initial vacancy
at a site $x$, then at any given
later time $t>0$ it is very likely to find  a vacancy in
$\bbZ^d_{x,\downarrow}\cup \{x\}$ close to $x$. These results will be the
keystone for the main outcome of this section (cf. Theorem \ref{th:1}), namely the fact that an
initial vacancy is able to generate a wave of equilibrium (\ie the
reversible measure $\pi$) in front of itself.  Finally, in Theorem \ref{th:2} we will estimate  the tail of the time needed to create a vacancy at  a given site for  the East--like process in the quadrant $\bbZ_+^d$. 
\subsection{Persistence of the vacancies}
Given $x\in \bbZ_+^d$ we say that $x$ is of class $n\in \bbN=\{0,1,\dots\}$ and
write $x\in \cC_n$
if $\min_{i}x_i\ge n$. Clearly $\cC_n\subset \cC_{n-1}$. 
\begin{lemma}
  \label{lem:0bis}
Let
$
p_n\triangleq \sup_{x\in \cC_n}\sup_{t\ge
  0}\sup_{\eta:\,\eta(x)=0}\bbP_{\eta}(\eta_t\restriction_{\L_x}=1),
$
where $\L_x=\prod_{i=1}^d[0,x_i]$. Then $p_n\le p^{n+1}$.
\end{lemma}
\begin{proof}
Let us fix
$x\in \cC_n$, $t>0$ and an initial  configuration $\eta$ such
that $\eta(x)=0$ and recall Definition \ref{def:1}. If we make the initial vacancy at $x$
``distinguished'', then $\{\eta_t\restriction_{\L_x}=1\}\subset  \{\cN_t\ge n+1\}\cap\{\eta_t(z^{(k)})=1\
\forall k=0,1,\dots,n\}$. Thus, using Proposition
\ref{prop:miracoli},
\[
\bbP_{\eta}(\eta_t\restriction_{\L_x}=1)\le p^{n+1}\bbP_\eta(\cN_t\ge
n+1)\le p^{n+1}.
\]
\end{proof}
The next is a simple but useful consequence of the above result. Fix
$\ell\ge 1$ and let $\cG_t$ be the event that there exists a
vertex $x\in V(\ell)\triangleq[-\ell+1,0]^d$ such that
$\cT_t(x)\ge t/\ell^d$, where
    \begin{equation}
      \label{eq:1}
\cT_t(x)\triangleq \int_0^t ds\, \mathbbold{1}{(\eta_s(x)=0)}      
    \end{equation}
is the total time 
spent in the zero state by the spin at $x$ up to time $t$.
    \begin{corollary} 
\label{cor:1bis} There exist positive constants $C,c$ such that
  \[    \sup_{\eta:\ \eta(0)=0}\bbP_\eta(\cG_t)\ge 1- C t\,\ell^d
  e^{-c \ell}, \quad \forall
\ell\ge 1.
\]
    \end{corollary}
    \begin{proof}
Observe that, if the box $V(\ell)$ was never completely filled during
the time-lag $t$, then necessarily
the event $\cG_t$ occurred. Thus it is sufficient to prove  that 
\[
\sup_{\eta:\ \eta(0)=0}\bbP_\eta(\exists\, s\le t:\
\eta_s\restriction_{V(\ell)}=1)\le C t\ell^de^{-c \ell},
\] 
for some constants $C,c>0$. Furthermore, 
using a union bound over the possible rings in
$V(\ell)$ within time $t$ (cf. the discussion after equation (5.14) in
\cite{East-Rassegna}) it suffices to prove that 
\[
    \sup_{\eta:\ \eta(0)=0}\sup_{s>0}\bbP_\eta(\eta_s\restriction_{V(\ell)}=1)\le e^{-c\ell}
\]
for some $c>0$. Such a bound follows from Lemma \ref{lem:0bis} with $c=-\log(p)$.
\end{proof}
\subsection{Local exponential ergodicity in $\bbZ^d$} 
\begin{theorem}
\label{th:1} There exist two positive constants $C,c$ such that the
  following holds. Fix $t>0$ and let $f$ be a function depending
  only on the spins in $[1,t^{1/2d}]^d$ with $\|f\|_\infty=1$.  Then
\[
\sup_{\eta:\ \eta(0)=0}|\bbE_\eta(f(\eta_t))-\pi(f)|\le C e^{-c t^{1/2d}}.
\]
\end{theorem}
\begin{corollary}
\label{cor:1}
Let $\nu$ be a probability measure on $\O=\{0,1\}^{\bbZ^d}$ such that,
 for all $\ell\ge 1$ and all $x\in \bbZ^d$, 
\[
\nu\left(\{\eta(y)=1\  \forall y\in [x,x+\ell-1]^d\}\right)\le e^{-m
  \ell},
\]
for some positive $m$. Then there exist $\l=\l(m)>0$ and $C>0$ such that
\[
\sup_x|\bbE_\nu\bigl(\eta_t(x)\bigr)-p| \le Ce^{- \l t^{1/2d}}.
\]    
\end{corollary}
\begin{proof}[Proof of the Corollary]
Without loss of generality consider only the case $x=0$ and, for
$\ell\ge 1$, write
\[
  |\bbE_\nu\left(\eta_t(0)\right)-p|\le 2 e^{-m \ell}+
  \int_{\{\eta:\ \exists x\in [-\ell,-1]^d:\ \eta(x)=0\}}
  d\nu(\eta)|\bbE_\eta\bigl(\eta_t(0)\bigr)-p|.
\]
Using Theorem \ref{th:1} the second term in the r.h.s. above is
smaller than $C e^{-c t^{1/2d}}$ for any $\ell \le t^{1/2d}$. Hence
the thesis.
\end{proof}
\begin{remark}
A similar result holds if one replaces the spin at
$x$ with an arbitrary function $f$ depending on finitely many
spins. In that case the constant $C$ will depend on $f$ through
$\|f\|_\infty$ and the size of the support of $f$ while the constant $\l$
will stay the same.    
\end{remark}
\begin{proof}[Proof of Theorem \ref{th:1}] In what follows $c$ will
  denote a generic constant depending on $p$ which may vary from line
  to line.  
Fix $f$ as in the theorem and assume for simplicity that
$\pi(f)=0$. Given $\eta$ such that
  $\eta(0)=0$, we use Corollary \ref{cor:1bis} with
  $\ell\equiv \lceil \d
  t^{1/2d}\rceil$ for some small positive constant $\d$ together with $\|f\|_\infty=1$ to write
\[
\bbE_\eta(f(\eta_t))=\bbE_\eta(f(\eta_t)\mathbbold{1}{(\cG_t)})+ O(e^{-c\d
  t^{1/2d}}),
\] 
where $\cG_t$ is the event that there exists $x\in [-\ell+1,0]^d$ with
$\cT_t(x)\ge t/\ell^d$. 
On the event $\cG_t$ let $\xi=(\xi_1,\dots,\xi_d)$ be the position of the smallest (in the
lexicographical order) vertex of
$V(\ell)=[-\ell+1,0]^d$ with the property that $\cT_t(\xi)\ge
t/\ell^d$ and let $\L_\xi$ be the box $[\xi_1+1,t^{1/2d}]\times
\prod_{i=2}^d[\xi_i,t^{1/2d}]$ (see Figure \ref{fig:1}).
\begin{figure}[t!]
\vspace{-0.275cm}
\includegraphics
[width=.45\textwidth]{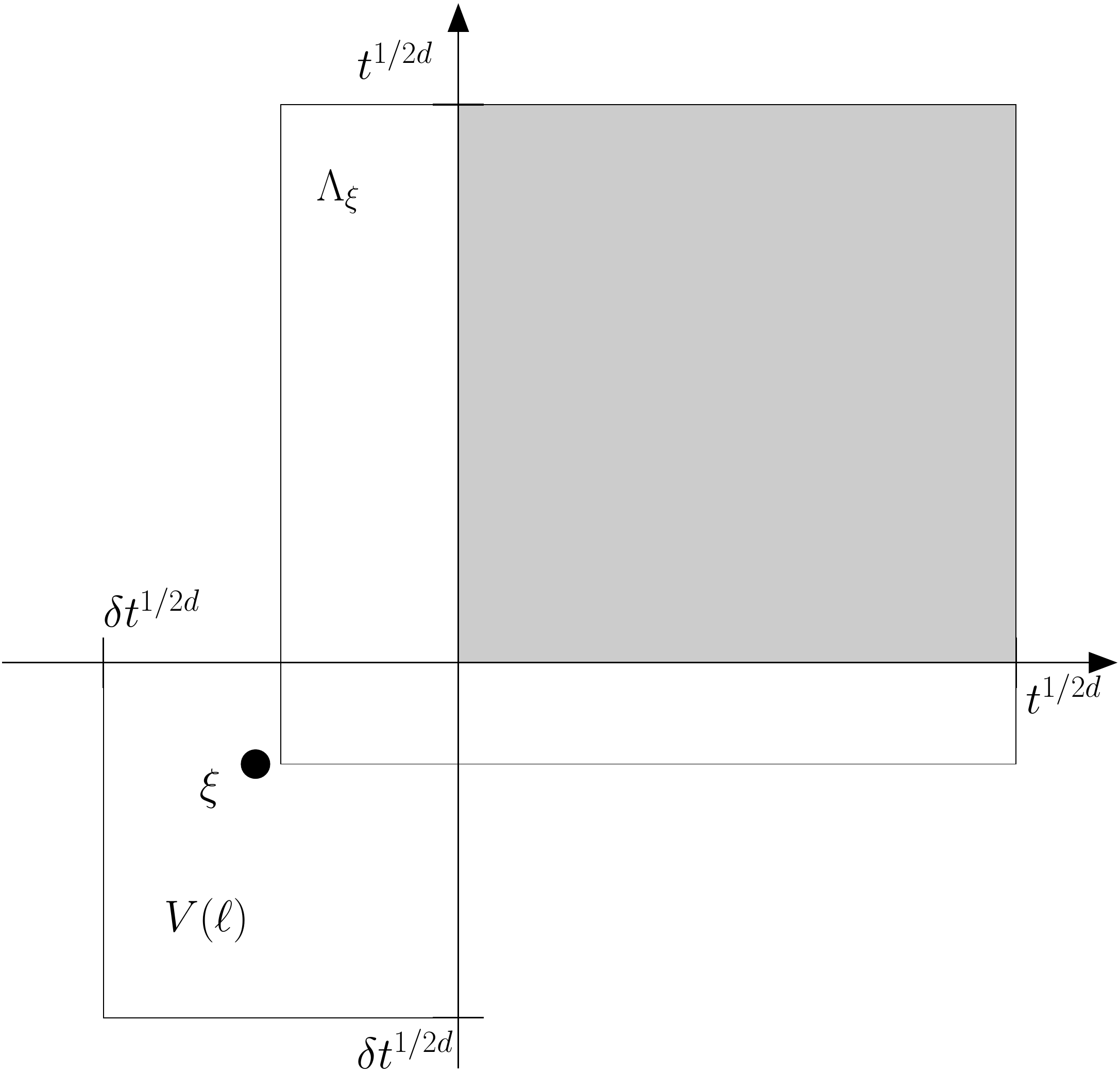}
\vspace{-0.25cm}
\caption{The box $\L_\xi$}
\label{fig:1}
\vspace{-0.5cm}
\end{figure}
 We then write
\[
\bbE_\eta(f(\eta_t)\mathbbold{1}{(\cG_t)})=\sum_{y\in
  V(\ell)}\bbE_\eta\left(\mathbbold{1}{(\hat \cG_{t,y})}\bbE_\eta(f(\eta_t)\tc
  \cF_{\L_y^c,t})\right),
\]
where $\hat \cG_{t,y}\triangleq\cG_t\cap \{\xi=y\}$. 
Notice that, given $\eta$, the event $\hat \cG_{t,y}$ is measurable w.r.t. $\cF_{\L_y^c,t}$.
 Clearly, for each $y\in V(\ell)$, one has
\begin{align*}
&|\bbE_\eta\left(\mathbbold{1}{(\hat \cG_{t,y})}\bbE_\eta(f(\eta_t)\tc \cF_{\L_y^c,t})\right)|\\
&\le \frac{1}{\min_{\s\in
    \O_{\L_y}}\pi(\s)}\bbE_\eta\Bigl(\mathbbold{1}{(\hat \cG_{t,y})}\sum_{\s
      \in \O_{\L_y}}\pi(\s)
      |\bbE_{\s\cdot \eta}(f(\s_t)\mid \cF_{\L_y^c,t})|\Bigr),
  \end{align*}
where $\s\cdot \eta$ means the configuration equal to $\s$ in $\L_y$
and to $\eta$ outside it. Above $\s_t\triangleq(\s\cdot \eta)_t
\restriction_{\L_y}$ and we used the fact that $f$ depends only
on the spins in $\L_y$.
We now observe that, given $\cF_{\L_y^c,t}$, the
evolution in $\L_y$ up to time $t$ is the standard East-like process with boundary
conditions that vary at times say $t_1,t_2,\dots,t_n$. Moreover, once the initial $\eta$ is given,
the times $\{t_i\}_{i=1}^n$ and the
actual value of the boundary conditions
$\{\eta_{t_i}(z)\}_{z\in \partial_E\L_y}$ become measurable
w.r.t. $\cF_{\L_y^c}$. 
Call $\cL^{(i)}$, $i=0,\dots,n$, the generator of the East-like process in $\L$ with
  boundary conditions given by $\eta_{t_{i}}$ (we set $t_0\equiv 0$ and $t_{n+1}\equiv t$). If $\eta_{t_{i}}(y)=0$,
then $\eta_{t_{i}} $ is an ergodic boundary condition and the
spectral gap of $\cL^{(i)}$ is not smaller than the spectral gap $\l$ in the
positive quadrant with minimal boundary conditions \cite{CFM2}.   Thus, using the Cauchy-Schwarz
inequality,
\begin{align}
\label{eq:100}
\sum_{\s
      \in \O_{\L_y}}\pi(\s)
      |\bbE_{\s\cdot \eta}(f(\s_t)\mid \cF_{\L_y^c,t})|&=\sum_{\s
      \in\O_{\L_y}}\pi(\s)\big|e^{t_1\cL^{(0)}}e^{(t_2-t_1)\cL^{(1)}}\dots e^{(t-t_n)\cL^{(n)}}f(\s)\big| 
\nonumber\\
&\le \|e^{t_1\cL^{(0)}}e^{(t_2-t_1)\cL^{(1)}}\dots
  e^{(t-t_n)\cL^{(n)}}f\|_\pi\nonumber \\
&\le  
 \exp\bigl[-\l
      \sum_{i=0}^n (t_{i+1}-t_i)\mathbbold{1}{(\eta_{t_i}(y)=0)}\bigr].
\end{align}
Above $\| \cdot \|_\pi$ denotes the norm in $L^2( \L_y, \pi)$. Note that 
in the last step we applied the classical inequality $\var(e^{t\cL}f)\le
e^{-2\gap(\cL)}\var(f)$, valid for any reversible continuous time
Markov chain, to the chains with generators $\{\cL^{(i)}\}_{i=0}^n$, 
together with  the fact  that $\pi\left(f^{(i)}\right)=0$ for all 
functions
\[
f^{(i)}\triangleq e^{(t_i-t_{i-1})\cL^{(i)}}\dots e^{(t-t_n)\cL^{(n)}}f\,.
\] 
This last property can be proved by induction on $i$ from $n$ to $0$ using that $\pi(f)=0$ and that $\cL^{(i)}$ is reversible (hence stationary) in $L^2(\L_y,\pi)$.

We now observe that $\sum_{i=0}^n (t_{i+1}-t_i)\mathbbold{1}{(\eta_{t_i}(y)=0)}=\cT_t(y)$, the total time spent in state zero by the spin at $y$. 
 By construction
the latter is at least $t/\ell^d\ge  \d^{-d}\sqrt{t}/2$ for $t$ large enough. In conclusion
\begin{gather*}
\frac{1}{\min_{\s\in
    \O_{\L_y}}\pi(\s)}\bbE_\eta\Bigl(\mathbbold{1}{(\hat \cG_{t,y})}\sum_{\s
      \in \O_{\L_y}}\pi(\s)
      |\bbE_{\s\cdot \eta}(f(\s_t)\mid \cF_{\L_y^c,t})|\Bigr)\\
\le 
\Bigl(\frac{1}{p\wedge q}\Bigr)^{|\L_y|}e^{-\l \d^{-d}\sqrt{t}/2}=O(e^{-c \sqrt{t}})
\end{gather*}
for $\d$ small enough and $t$ large enough. Thus 
\[
\big|\bbE_\eta(f(\eta_t)\mathbbold{1}{(\cG_t)})\big|\le C\ell^d \,e^{-c \sqrt{t}}
\]
for some constant $C$ and the result follows.
\end{proof}
\subsection{Exponential tail of the persistence times in $\bbZ^d_+$}
\label{sec:persistence-time}
Consider the East-like process in $\bbZ^d_+$ with minimal boundary
condition, and let $\bbP_\eta^{\rm min}(\cdot)$ denote its law when
the initial configuration is $\eta$. Recall also definition
\eqref{eq:1} of the random variable $\cT_t(x)$.
\begin{theorem} 
\label{th:2} 
 There exist $\kappa,\l>0$ and $\d\in (0,1)$  such that the following
holds. For all $x\in \bbZ^d_+$ and $t\ge \kappa \|x\|_1$,
\begin{equation}
  \label{eq:2}
  \sup_\eta\bbP^{\rm min}_{\eta}(\cT_t(x) \le \d^d t) \le d e^{-\l \d  t}.
\end{equation}
In particular
\begin{equation}
  \label{eq:3}
\sup_\eta\bbP^{\rm min}_{\eta}(\t_x \ge t) \le de^{-\l\d  t},
\end{equation}
where $\t_x\triangleq \inf\{t:\ \eta_t(x)=0\}$.
\end{theorem}
\begin{remark}
For simplicity we have stated  the result with minimal boundary
conditions. Actually the
same proof, with minor modifications, holds for any
ergodic boundary conditions with uniform constants $\kappa,\l,\d$.
\end{remark}
\begin{proof}
It is obvious that \eqref{eq:3} follows from \eqref{eq:2}. To prove
the latter 
we first need a technical lemma. Let $x\in \bbZ^{d}_+$ and write $x=(x^*,x_d)$ where $x^*=(x_1,\dots,x_{d-1})\in \bbZ_+^{d-1}$. Let $I=\{y=(x^*,j),\
j=1,\dots x_d\}$ if $x_d>0$ and $I=\emptyset$ otherwise (see Figure \ref{fig:2}).
\begin{figure}[!ht]
\begin{overpic}[scale=0.3]{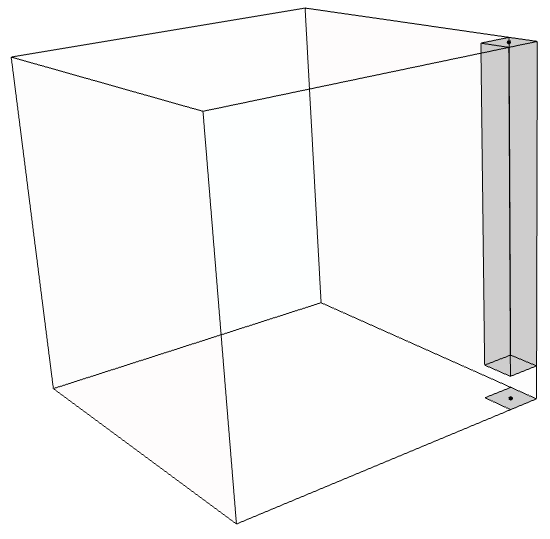} 
\put(-3,18){$(0,0)$}
\put(90,95){$x=(x_1,x_2,x_3)$}
\put(95,18){$(x^*,0)$}
\end{overpic}
\caption{The interval $I$ and the point $x^*$ in three dimensions.}
\label{fig:2}\end{figure}
 For the East-like process in $\bbZ^d_+$
with ergodic boundary conditions and initial configuration $\eta$, the variable $\cT_t((x^*,0))$ is
measurable w.r.t. $\cF_{I^c,t}$.
Moreover, since the projection on the slab $\bbZ_+^{d-1}\times
\{0\}$ of
the East-like process in $\bbZ^d_+$ with minimal boundary conditions
coincides with the East-like process on $\bbZ^{d-1}_+$ with
minimal boundary conditions, the law of $\cT_t((x^*,0))$ under
$\bbP_\eta^{\rm min}$ coincides with the law of $\cT_t(x^*)$ in
$(d-1)$ dimensions. In the sequel, for notation convenience, we will
simply write $\cT_t(x^*)$ instead of the more precise $\cT_t((x^*,0))$.
\begin{lemma}
\label{lem:2}
There exist constants $\d\in (0,1)$ and $c>0$ such that 
\[
\sup_{\eta}\bbP_\eta^{\rm min}(\cT_t(x)\le \d\cT_t(x^*)\tc \cF_{I^c,t})\le (p\wedge
  q)^{-x_d}\  e^{-c \cT_t(x^*)}.
\]  
\end{lemma}
We postpone the proof of the lemma and write
\begin{align*}
 \sup_{\eta}\bbP^{\rm min}_\eta\left(\cT_t(x)\le \d^d t\right)
\le  &\sup_{\eta}\bbP^{\rm min}_\eta\left(\cT_t(x^*)\le
  \d^{d-1} t\right) \\&+
 \sup_{\eta}\bbP^{\rm min}_\eta\left(\cT_t(x)\le \d^d t\,;\,\cT_t(x^*)\ge \d^{d-1} t \right)\\
&\le  \sup_{\eta}\bbP^{\rm min}_\eta\left(\cT_t(x^*)\le \d^{d-1}t\right)+ (p\wedge
  q)^{-x_d}\  \ e^{-c \d^{d-1}t},
\end{align*}
where we used Lemma \ref{lem:2} to bound the second term in the
r.h.s. of the first inequality. Notice that the first term in the
r.h.s above has the
same form of the starting quantity but now in $(d-1)$ dimensions. We can then iterate on the dimension
$d$ to get
\[
  \sup_{\eta}\bbP^{\rm min}_\eta\left(\cT_t(x)\le \d^d t\right)\le 
 \sup_{\eta}\bbP^{\rm East}_\eta\left(\cT_t(x_1)\le \d t\right) +\sum_{i=2}^{d}(p\wedge
  q)^{-x_{i}}\  e^{-c \d^{i-1}t},
\]
where $\bbP^{\rm East}_\eta(\cdot)$ is the law of the one dimensional
East process on the interval $[0,x_1]$ with ergodic boundary
conditions. As in \cite{CMRT}*{Theorem 3.6}, the
first term in the r.h.s. above can be bounded by 
\[
 \sup_{\eta}\bbP^{\rm East}_\eta(\cT_t(x_1)\le \d t)\le (p\wedge
 q)^{-x_1}\bbP^{\rm East}_\pi(\cT_t(x_1)\le \d t)\le (p\wedge
 q)^{-x_1}\, e^{-c' \d t},
\]  
for some constant $c'>0$ provided that $\d$ is small enough
(independent of $t$). The proof of \eqref{eq:2} is finished by choosing e.g. $\kappa =
2\log\left(1/(p\wedge q)\right)/(c\wedge c') \d $ and $\l=(c\wedge c')/2$.
\end{proof}

\begin{proof}[Proof of Lemma \ref{lem:2}]
Using the exponential Chebyshev inequality
we get
\begin{align*}
\sup_{\eta}\bbP_\eta^{\rm min}\Bigl(\cT_t(x)&\le \d\cT_t(x^*)\tc \cF_{I^c,t}\Bigr)\\
&\le \sup_\eta \inf_{\g>0}e^{-\g(t-\d\cT_t(x^*))}\  
\bbE_\eta^{\rm min}\left(e^{\g(t-\cT_t(x))}\tc \cF_{I^c,t}\right)\\
&\le (1/p\wedge q)^{x_d}\sup_\eta
  \inf_{\g>0}e^{-\g(t-\d\cT_t(x^*))}\ 
\bbE_{\eta,\pi}^{\rm min}\left(e^{\g(t-\cT_t(x))}\tc\cF_{I^c,t}\right),
\end{align*}
where $\bbE_{\eta,\pi}^{\rm min}(\cdot)$ denotes the expectation
w.r.t. the East-like process with initial law
$\nu(\eta')=\pi(\eta'\restriction_{ I})\mathbbold{1}{(\eta' \restriction_{ \bbZ^d_+\setminus
  I}=  \eta \restriction_{ \bbZ^d_+\setminus
  I}  )}$. 

We will now bound the term 
$
\bbE_{\eta,\pi}^{\rm min}\left(e^{\g(t-\cT_t(x))}\tc
\cF_{I^c,t}\right)=\bbE_{\eta,\pi}^{\rm min}\left(e^{\g\int_0^t ds \, \mathbbold{1}{(\eta_s(x)=1)}}\tc
\cF_{I^c,t}\right)
$ 
using the Feynman-Kac formula. Firstly notice that it is enough to compute
the above expectation w.r.t. the process in the box
$\L_x=\prod_{i=1}^d[0,x_i]$. Secondly, since the East boundary of any
vertex in $\L_x\setminus I$ does not intersect $I$ and since the
initial configuration in $\L_x\setminus I$ is deterministic, the
projection of the process in $\L_x\setminus I$ is measurable
w.r.t. $\cF_{I^c,t}$.  
Denote by
$0<t_1<t_2,\dots, <t_n< t$ the successive times in $[0,t]$ at which one of the spins at the
East boundary of the interval $I$ changes. During any time interval of the
form $[t_i,t_{i+1})$ (define $t_0\equiv 0$ and $t_{n+1}\equiv t$) the
process in $I$ is nothing but the usual one dimensional East
process with possibly certain vertices which are unconstrained, namely
those with a zero spin at their East boundary (in particular,  the initial law
$\pi$ in $I$ is preserved at any later time). Let $\cL^{(i)}$ be the
corresponding Markov generator.  
 Because of what we just said
$\cL^{(i)}$ is self-adjoint in $L^2(I,\pi)$ and, if the boundary spin
at $(x^*,0)$ is zero, it  has also a spectral gap which is not smaller
than the spectral gap  $\gap(\cL_{\rm East})$ on the East process on $\bbZ$ \cite{CMRT}. Moreover, for any function
$F:\O_I\mapsto \bbR$ and any $s>0$, the Feynman-Kac formula
\[
\bbE^{(i)}_\s\left(F(\s_{t_{i+1}-t_i})e^{\g \int_0^{t_{i+1}-t_i} ds\,
    \mathbbold{1}{(\eta_{s}(x)=1)}}\right)= \left(e^{(t_{i+1}-t_i)(\cL^{(i)}+\g \cA_x)}F\right)(\s)
\]   
holds, with $\cA_xF (\s)= \mathbbold{1}{(\sigma(x)=1)}F(\s)$ and
$\bbE^{(i)}_\s(\cdot)$ being the expectation over the process with
generator $\cL^{(i)}$ and initial condition $\s$. 

Thus
\begin{align}
  \label{eq:4}
  \bbE_{\eta,\pi}^{\rm min}&\left(e^{\g\int_0^t ds \,
      \mathbbold{1}{(\eta_s(x)=1)}}\tc
\cF_{I^c,t}\right)\nonumber \\
&= \langle {\mathbf 1},e^{t_1(\cL^{(0)}+\g
    \cA_x)}\cdot e^{(t_2-t_1)(\cL^{(1)}+\g \cA_x)}\cdot\dots \cdot
                                 e^{(t-t_n)(\cL^{(n)}+\g
                                 \cA_x)}{\mathbf 1}\rangle_\pi\nonumber\\
&\le \prod_{i=0}^n \|e^{(t_{i+1}-t_i)(\cL^{(i)}+\g \cA_x)}\|_\pi
\end{align}
where $\langle\cdot\rangle_\pi$ denotes the scalar product in
$L^2(I,\pi)$ and $\|\cdot\|_\pi$ denotes the operator norm on $L^2(I,\pi)$. 

For $i$ such that $\eta_s((x^*,0))=1,\ s\in [t_i,t_{i+1})$,
we simply bound $\|e^{(t_{i+1}-t_i)(\cL^{(1)}+\g \cA_x)}\|_\pi$ by $e^{\g
  (t_{i+1}-t_i)}$. Indeed, $\cL^{(1)}+\g\cA_x\le \g\mathds 1$ as self-adjoint
operators. In the opposite case and provided that $\g$ is
small enough (e.g. $\g=\gap(\cL_{\rm East})/2$) we can use a result 
from \cite{CMRT}*{Proof of Theorem 3.6} to get that 
\[
 \|e^{(t_{i+1}-t_i)(\cL^{(1)}+\g \cA_x)}\|_\pi\le e^{ \b \g(t_{i+1}-t_i)}
\]
where $\b=pq/(1+p)+p<1$. 

In conclusion
\[
 \prod_{i=0}^n \|e^{(t_{i+1}-t_i)(\cL^{(1)}+\g \cA_x)}\|_\pi \le e^{\b \g
   \cT_t(x^*) +\g (t-\cT_t(x^*))}
\]
and 
\[
\sup_{\eta}\bbP_\eta^{\rm min}\Bigl(\cT_t(x)\le \d\cT_t(x^*)\tc
\cF_{I^c,t}\Bigr) \le  (1/p\wedge q)^{x_d} e^{-\g(1-\b-\d)\cT_t(x^*)}.
\]
The proof is complete if we take e.g. $\d= (1-\b)/2$.
\end{proof}

\section{Mixing time of the East-like process}
Consider the East-like process in $\L=[1,L]^d$ with ergodic boundary
conditions and let 
\[
\tmix\triangleq \inf\{t>0:\ \max_\eta \|\bbP^\L_\eta(\eta_t\in
\cdot)-\pi\|_{\rm TV}\le 1/4\}
\]
be its mixing time. 
\begin{theorem}
\label{th:mix}
There exists $C>0$ such that $C^{-1}L\le \tmix\le CL$ for all $L \geq 1$.
\end{theorem}
\begin{remark}
  Notice that the standard inequality $\tmix
  \le \text{const. }\log(1/\pi^*)\times \trel $ with $\pi*=\min_{\eta\in \O_\L}\pi(\eta)$   
(see e.g. \cite{Saloff}) only gives $\tmix =O(L^d)$. Also, if the boundary
conditions are minimal, it is easy
to see that the logarithmic
Sobolev constant $\a_\L$ (cf. \cite{Saloff}) satisfies $\a_\L\sim L^{-d}$. The upper
bound is proved by plugging into the variational
characterization of $\a_\L$ the test function given by the indicator of the configuration without
vacancies. The lower bound follows at once from the general bound
$\a_\L\ge
\gap(\cL_\L^{\rm min})/(2+\log(1/\pi^*))$.   
Thus the East-like process with
minimal boundary conditions has $\trel=O(1)$, $\tmix \sim L$ and a
logarithmic Sobolev constant $\sim L^{-d}$.    
\end{remark}
 We first establish three key preliminary results before proving Theorem \ref{th:mix}. We denote by
$\hat \O_\L$ the set of configurations such that in \emph{any} interval 
$I\subset \L$ parallel to one of the coordinate axes and of length
$\lfloor(\log L)^2\rfloor $
there exists at least one vacancy. The
first result says that any initial
configuration will, with high probability, evolve into $\hat \O_\L$
in a time $t=O(L)$. 
\begin{lemma}
\label{lem:1}
For any $\e\in (0,1)$ there exists $M$ such that, for all initial
configurations $\eta\in \O_\L$ and all $L\geq 1$,
\[
\sup_{t\ge ML}\bbP^\L_\eta \left(\eta_t\notin
  \hat \O_\L\right)\le \e. 
\]   
\end{lemma}
\begin{proof} Fix an interval $I\subset [1,L]^d$ of the form
  $I=[x,x+\lfloor (\log L)^2\rfloor e]$, $e\in \cB$. 
Using Theorem \ref{th:2} and its notation, there exist $M>0$ and
$c>0$ independent of $L$ and of
the initial configuration such
that
\[
\bbP^\L_\eta(\tau_x\ge ML/2)\le e^{-c ML/2}.
\]
The strong Markov property at the hitting time $\t_x$ together
with Theorem \ref{th:1} gives that, for $t\ge ML $, 
\begin{gather*}
\bbP^\L_\eta(\t_x\le
  ML/2\,;\, \eta_t(z)=1\ \forall \, z\in I)\\
=
\bbE^\L_\eta\left(\mathbbold{1}{(\t_x\le
  ML/2)}\bbP^\L_{\eta_{\t_x}}(\eta_{t-\t_x}(z)=1\ \forall \, z\in I) \right)\\
\le \pi(\eta(z) =1\ \forall \, z\in I) + Ce^{-c (t-ML/2)^{1/2d}}= O(e^{-c'
  (\log L)^2}).  
\end{gather*}
Since the number of such intervals $I$ is $O(e^{c'' \log L})$, a union
bound over the choice of $I$ finishes the proof.
\end{proof}
The next result is a small refinement of the arguments used in the proof of
Theorem \ref{th:1}. 
Consider the East-like process in $\bbZ^d_+$ with ergodic
  boundary conditions and let  $\cG_{t}(x;\D)$ be
  the event that the spin at $x$ was unconstrained (\ie with at least
  a vacancy in its East-boundary) for a total time
  $0\le \D\le t$ during the time interval $[0,t]$. By
  construction, given the initial $\eta$,   $\cG_{t}(x;\D)$ is measurable w.r.t.  $\cF_{\bbZ^d_{x,\downarrow},t}$.
\begin{lemma} 
\label{lem:3}
Given $x\in \bbZ_+^d$, let $V$ be a box of
  the form 
$V=\prod_{i=1}^d
  [x_i,x_i+\ell-1]$,  $\ell\geq 1$,  
and let $f:\bbZ^d_+\mapsto \bbR$ be a bounded function which
does not depend on the spins in $\bbZ^d_{x,\uparrow}\setminus V$. There
  exist positive constants $c,\l$ so that,
  on the event  $\cG_{t}(x;\D)$,
\[
\max_\eta|\bbE^{\bbZ^d_+}_{\eta}\bigl(f(\eta_t)-f^V(\eta_t)\,\tc
\cF_{V^c,t}\bigr)| \le C \|f\|_\infty \,e^{c\ell^d- \l \D},
\]
where $f^V(\eta)\triangleq \pi_V(f)(\eta)$ is the equilibrium average
in $V$ of $f$.
The constant $\l$ can be chosen as the spectral gap of the process in
$\bbZ^d_+$ with
minimal boundary conditions.  
\end{lemma}
\begin{proof}
Let $W=V\cup (\bbZ^d_+\setminus \bbZ^d_{x,\uparrow})$ (see Figure
\ref{fig:4}). 
\begin{figure}[!ht]
\vspace{-0.275cm}
\begin{overpic}[scale=0.3]{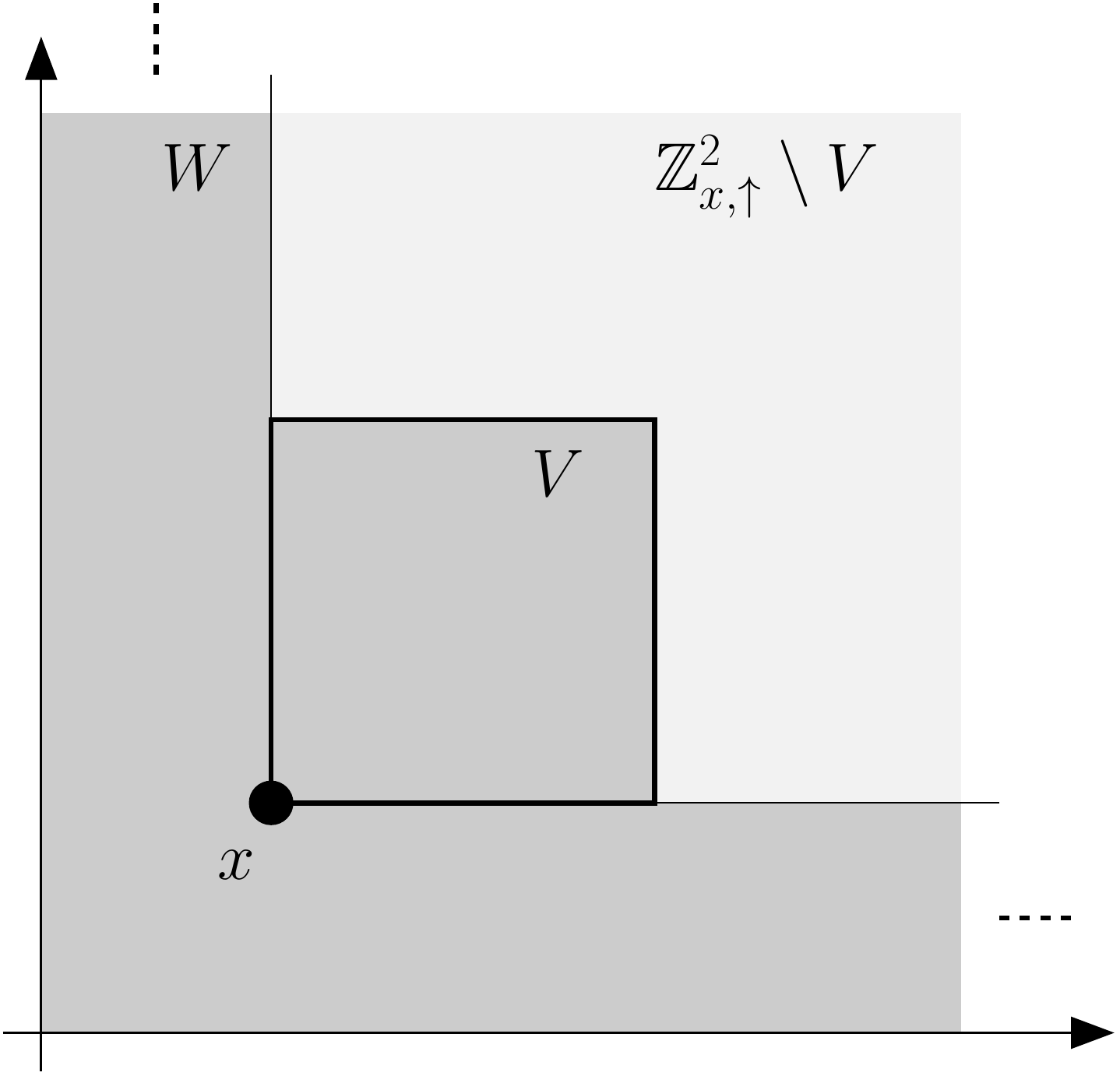}  
\end{overpic}
\caption{The set $W$ when $d=2$.}
\label{fig:4}
\end{figure}

Using the assumption that $f$
does not depend on the spins in $\bbZ^d_{x,\uparrow}\setminus V$ and the
oriented nature of the constraints,  
we can safely replace the East-like process in $\bbZ^d_{+}$ with the East-like
process in $W$, i.e.
\[
\bbE^{\bbZ^d_{+}}_\eta\bigl(f(\eta_t)-f^V(\eta_t)\tc
\cF_{V^c,t}\bigr)=\bbE^{W}_{\eta}\bigl(f(\eta_t) -f^V(\eta_t)\tc
\cF_{V^c,t}\bigr).
\] 
Clearly, given the initial $\eta$,  the dynamics of
the spins in $W\setminus V$ is measurable
w.r.t. $\cF_{V^c,t}$. Moreover, using Proposition \ref{lem:0}, 
\[
\pi_V\left(\bbE^{W}_{\eta}\left(f(\eta_s)\tc
    \cF_{V^c,s}\right)\right)=f^V(\eta_s) \quad a.s.\quad \forall\
s\in [0,t],
\]
where $\pi_V$ averages only over the spins in $V$.
The proof now follows the pattern of the proof of
Theorem \ref{th:1} (cf. \eqref{eq:100}). We write
\begin{align*}
&\mathbbold{1}{(\cG_{t}(x;\D))}
  |\bbE^W_\eta\left(f(\eta_t)-f^V(\eta_t)\tc \cF_{V^c,t}\right)| \\
&\le
1/(p\wedge q)^{\ell^d}
\mathbbold{1}{(\cG_{t}(x;\D))}\pi_V\left(|\bbE^W_\eta\left(f(\eta_t)-f^V(\eta_t)\,\tc
    \cF_{V^c,t}\right)|\right)\\
&\le
1/(p\wedge q)^{\ell^d}
\mathbbold{1}{(\cG_{t}(x;\D))}\bigl[\var_{\pi_V}\bigl(\bbE^W_\eta\left(f(\eta_t)\tc
    \cF_{V^c,t}\right)\bigr)\bigr]^{1/2}
\le \|f\|_\infty\, e^{c\ell^d - \l\D},
\end{align*}
where $\var_{\pi_V}$ denotes the variance w.r.t to $\pi_V$ and
$c=-\log(p\wedge q)$.
The fact that $\l$ can be taken equal to the spectral gap of the first quadrant with
minimal boundary conditions follows immediately from the monotonicity
of the spectral gap in the boundary conditions and in the volume.  \end{proof}

The third result says that the East-like process in the box
$\L=[1,L]^d$ with ergodic boundary conditions reaches equilibrium (in total variation distance)
in a time lag $O(\log(L)^{4d})$ if the initial
configuration belongs to  $\hat \O_\L$. More precisely, let $\mu^\eta_t$ denote the
law of $\eta_t$ under $\bbP^\L_\eta(\cdot)$. 
\begin{lemma}
\label{lem:4}
For any $\e\in (0,1)$ there exists $L_0$ such
that the following holds. Let $T=(\log L)^{5d}$ and let $
d(T)=\max_{\eta\in \hat \O_{\L}}\|\mu^\eta_T-\pi_{\L}\|_{\rm TV}.
$
Then 
$
\sup_{L\ge L_0}d(T)\le \e.
$   
\end{lemma}
\begin{proof}
Fix $\eta\in \hat \O_\L$ and let us order the vertices of $\L=[1,L]^d$ as
follows. We first choose some order of the vertices on each hyperplane
$\cP=\{x\in \bbZ^d_{+}: \ \|x\|_1= \text{const.}\}$ in such a way that two
consecutive vertices have distance $O(1)$ for large $L$ (e.g. they have exactly two coordinates where they differ by
one). Then, for any pair
$x,y\in \L$, we say that $x\prec
y$ if either $\|x\|_1<\|y\|_1$ or if $x$ comes before $y$ when they
have the
same $\ell^1$-norm. The $i^{th}$ vertex in the above order will be
denoted by $x^{(i)}$.

Next, for any
$f:\O_\L\mapsto \bbR$ with $\|f\|_\infty=1$ and any $j=1,\dots, n$,
where $n=L^d$, we denote
by $f^{(j)}$ the new function obtained by averaging $f$ w.r.t. the
equilibrium measure $\pi$ over the last
$j$ spins in the above order. For notation convenience we set
$f^{(0)}\equiv f$. Notice that $f^{(j)}=\pi_{j}(f^{(j-1)})$, where $\pi_j(\cdot)$ denotes
the marginal of $\pi$ over the $(L^d-j+1)^{th}$-spin, and that
$f^{(n)}=\pi(f)$. Then
\[
|\mu^\eta_t(f)-\pi(f)|\le
\sum_{j=1}^{n}|\mu^\eta_t (f^{(j-1)})-\mu^\eta_t (f^{(j)})|,\quad \forall t\ge 0
.
\]
We now choose $t=(\log L)^{5d}$ and $\eta\in
\hat \O_\L$ and prove that each term in the above sum is smaller
than $O(e^{-(\log L)^2})$ for large enough $L$. 

Let $\L_j$ be the set  $\{x^{(1)},x^{(2)},\dots,x^{(n-j)}\}$.
\begin{figure}[!ht]
\begin{overpic}[scale=0.3]{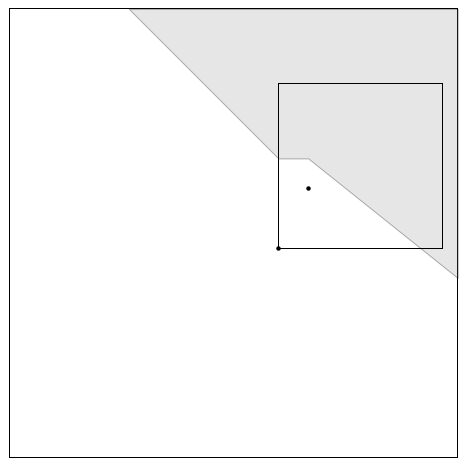}
\put(55,42){$z$}
\put(66,55){$x^{(n-j+1)}$}
\put(85,75){$V$}
\end{overpic}
\caption{The shaded quasi-triangle corresponds to vertices larger than
$x^{(n-j+1)}$}
\label{fig:3}\end{figure}
Using the assumption that
the initial configuration $\eta$ belongs to $\hat \O_\L$, there
exists a vertex $y\in (\L_j\cup \partial_{E}\L_j)\cap \bbZ^d_{x^{(n-j+1)},\downarrow}$ such that
$|y-x^{(n-j+1)}|\le (\log L)^2$ and $\eta(y)=0$.
Corollary \ref{cor:1bis} implies that, with probability greater than $1-e^{-(\log
  L)^2}$, there exists $C>0$ independent of $j$ and a vertex $z\in (\L_j\cup \partial_{E}\L_j)\cap \bbZ^d_{x^{(n-j+1)},\downarrow}$ such that:
\begin{enumerate}[(i)]
\item $|z-x^{(n-j+1)}|\le C(\log L)^2$,
\item one spin in the East boundary of $z$ spends a time greater than $(\log L)^{3d}$ in state $0$ up to 
  time $t$. 
\end{enumerate}
Notice that,  if $V\equiv \prod_{i=1}^d[z_i,z_i+3C(\log
L)^2]$, then neither $f^{(j)}$ nor $f^{(j-1)}$ depend on the spins in
$\bbZ^d_{z,\uparrow}\setminus V$ (see Figure \ref{fig:3}) and we can apply Lemma \ref{lem:3} to both $(V,f^{(j-1)})$
and $(V,f^{(j)})$. Using the fact that
$\pi_{V}(f^{(j-1)})=\pi_{V}(f^{(j)})$ we finally get that
\begin{align*}
  |\mu^\eta_t(f^{(j-1)})-\mu^\eta_t(f^{(j)})|\le e^{-(\log  L)^2}+ C'
  e^{c (\log L)^{2d}-\l(\log L)^{3d}}.  
\end{align*}
In conclusion, for any $L$ large enough
\[
|\mu_t(f)-\pi(f)|\le 2  L^d e^{-(\log L)^2}
\]
and the lemma follows.
\end{proof}
We are finally in a position to prove our main result.
\begin{proof}[Proof of Theorem \ref{th:mix}]
The lower bound is straightforward by choosing as initial condition
the configuration without vacancies and using the finite speed of
propagation of information (see e.g. \cite{East-cutoff}*{Section 2.4}) to prove that, with high
probability, the process is not able to create vacancies near the
vertex $\hat x=(L,\dots,L)$ in a time $t=\d L$, if $\d$ is small enough.  

To prove the upper bound we proceed as follows. Using Lemma \ref{lem:1}, uniformly in the initial condition, in a time $t=O(L)$ the East-like chain will
enter the good set $\hat \O_\L$ with probability e.g. greater than
$7/8$. Using the Markov property and Lemma \ref{lem:4}, in an additional  time lag
$O((\log L)^{5d})$ the chain will reduce its variation distance from
the target distribution $\pi_\L$ to less than $1/8$. 
\end{proof}

\section*{Acknowledgments}
F. Martinelli is grateful to the Centre for Mathematics and
Computer Science (CIMI) in Toulouse for their kind invitation to the
workshop `Talking Across Fields" and for their financial support. We
would also like to thank the organizers of the   TAU-CECAM
workshop on ``Percolation and the Glass Transition'' (Tel Aviv 2014) where part of these results were discussed. 
P.C. would like to acknowledge the support
of the University of Warwick IAS through a Global Research Fellowship.

\end{document}